\documentclass[final]{AAM}
\usepackage{xcolor}
\usepackage{mathrsfs}

\newcommand{\TheAuthor}{}
\newcommand{\Author}[1]{\renewcommand{\TheAuthor}{#1}}   

\newcommand{\TheTitle}{}
\newcommand{\Title}[1]{\renewcommand{\TheTitle}{#1}}     

\newcommand{\TheDateAccepted}{}
\newcommand{\DateAccepted}[1]{\renewcommand{\TheDateAccepted}{#1}}

\newcommand{\TheDateReceived}{}
\newcommand{\DateReceived}[1]{\renewcommand{\TheDateReceived}{#1}}

\usepackage{fancyhdr}                    
\renewcommand{\headrulewidth}{0.0pt}

\Author{xxx}
\Title{AAM: Intern. J., Vol. 1, Issue 1 (January 1900)}
\DateReceived{Aug. 18, 2021}
\DateAccepted{Sept. 30, 2021}

\begin{document}


\def\undertext#1{$\underline{\smash{\hbox{#1}}}$}

\renewcommand{\headrulewidth}{0pt} 

\title{A Schwartz-Zippel Type Estimate for Homogenous Finite Field Polynomials}
\maketitle

\centerline{\bf  $^1$Ghurumuruhan Ganesan}

\begin{center}
Institute of Mathematical Sciences, HBNI, Chennai,
Chennai, India\\
$^1${\color[rgb]{0,0,1}{\undertext{gganesan82@gmail.com}}} 
\end{center}

\vspace{0.5truecm}
\centerline{\received\TheDateReceived; \accepted\TheDateAccepted}

\section*{Abstract}
In this paper, we obtain a Schwartz-Zippel type estimate for homogenous finite field polynomials. Specifically, we use a probabilistic recursion technique to find upper and lower bounds for the number of zeros of a homogenous polynomial  and illustrate our result with two examples involving perfect matching in bipartite graphs and common zeros in a collection of polynomials, respectively.




{\Large{\bf{MSC 2010 No.: }}}{11T06}
\noindent

\renewcommand{\theequation}{\thesection.\arabic{equation}}
\setcounter{equation}{0}
\section{Introduction} \label{intro}

The Schwartz-Zippel bound~\cite{schwa} and combinatorial nullstellensatz~\cite{alon} are important tools from both theoretical and application perspectives. The Schwartz-Zippel bound provides an estimate on the number of zeros of a polynomial whose coefficients take values in a finite field and has applications in polynomial comparison, primality testing, perfect matching in graphs etc. The combinatorial nullstellensatz determines conditions under which a polynomial contains a non-zero in a given set whose cardinality is larger than the overall degree of the polynomial. We refer to Chapter 7,~\cite{mot} for more material. Recently~\cite{hossan} have used algebraic techniques to obtain Schwartz-Zippel type bounds for intersection of algeraic varieties with Cartesian products of two-dimensional sets.

In this paper, we use probabilistic methods and recursive techniques to find upper and lower bounds the number of zeros of a homogenous polynomial. As an illustration of the upper bound, we  apply our bounds to estimate the error probability in the randomized algorithm that determines the presence of a perfect matching in a bipartite graph. We use the lower bound on the number of zeros to determine the presence of common zeros in a set of polynomials.

\subsection*{Homogenous polynomials}
For integers~\( l\geq 1, q = p^{l}, p\) prime, let~\(\mathbb{F}_q\) be the finite field consisting of~\(q\) elements with characteristic~\(p.\) A homogenous polynomial in~\(\mathbb{F}_q\) in the variables~\(x_1,\ldots,x_m\) and having degree~\(k \leq m\) is of the form
\begin{equation}\label{poly}
Q(x_1,\ldots,x_m) = \sum_{I : \#I = k} \delta_I \cdot \prod_{i \in I} x_i
\end{equation}
where~\(\delta_I \in \mathbb{F}_q\)  are not all zero and the sum is over all subsets~\(I \subseteq \{1,2,\ldots,m\}\) of cardinality~\(\#I = k.\)

We say that the~\(m-\)tuple~\((y_1,\ldots,y_m) \in \mathbb{F}_q^{m}\) is a \emph{zero} for~\(Q\) if~\[Q(y_1,\ldots,y_m) = 0.\]
\begin{theorem}\label{thm1} If~\({\cal Z} \subseteq \mathbb{F}^{m}_q\) denotes the set of zeros of~\(Q,\) then
\begin{equation}\label{z_est}
q^{m-k} (q-1)^{k-1} \leq \#{\cal Z} \leq q^{m-k}(q^{k}-(q-1)^{k})
\end{equation}
where~\(\#{\cal Z}\) denotes the cardinality of the set~\({\cal Z}.\)
\end{theorem}
Using~\((1+t)^{k} - t^{k} = \sum_{l=0}^{k-1} {k \choose l} t^{l} \leq t^{k-1}(2^{k}-1)\) with~\(t = q-1,\) we get from~(\ref{z_est}) that
\[q^{m-k} (q-1)^{k-1} \leq \#{\cal Z} \leq q^{m-k}(q-1)^{k-1}(2^{k}-1).\] Thus the number of zeros is of the order of~\(q^{m-k}(q-1)^{k-1}.\)

We use recursion to prove Theorem~\ref{thm1} in the next Section. If~\(k\leq q,\) we use~\(\left(1-\frac{1}{q}\right)^{k} \geq 1-\frac{k}{q}\) to get from~(\ref{z_est}) that~\(\#{\cal Z} \leq k \cdot q^{m-1},\) the Schwartz-Zippel bound~\cite{mot}. We remark that in Theorem~\ref{thm1}, we can choose the field size~\(q\) independent of the degree~\(k\) of the polynomial and the number of variables~\(m.\) In Section~\ref{examples}, we describe applications for both the lower and the upper bounds.

The paper is organized as follows. In Section~\ref{pf_thm1}, we prove Theorem~\ref{thm1} and in Section~\ref{examples} we describe applications that use the bounds in~(\ref{z_est}).


\setcounter{equation}{0}
\renewcommand\theequation{\thesection.\arabic{equation}}
\section{Proof of Theorem~\ref{thm1}} \label{pf_thm1}
We use the probabilistic method analogous to the proof of Schwartz-Zippel Lemma and obtain a recursion relation
involving the probability that the set of randomly chosen values indeed form a zero of the polynomial.
We then evaluate the recursive relation to obtain the desired bounds on the corresponding probability.

We begin with the upper bound in~(\ref{z_est}). Let~\(r_1,\ldots,r_m\) be independently and uniformly randomly chosen from~\(\mathbb{F}_q\) and let
\begin{equation}\label{pk_def}
p_{k,m} := \max_{Q} \mathbb{P}(Q(r_1,\ldots,r_m) = 0)
\end{equation}
where the maximum is over all homogenous polynomials in~\(m\) variables and with degree~\(k.\)
To prove the upper bound in~(\ref{z_est}) it suffices to see that
\begin{equation}\label{pk_est_gen}
p_{k,m} \leq 1 - \left(1-\frac{1}{q}\right)^{k}.
\end{equation}

We first prove~(\ref{pk_est_gen}) for~\(k=1.\) Consider the polynomial~\[Q(x_1,\ldots,x_m) = \sum_{i=1}^{m} \alpha_i \cdot x_i,\]
where~\(\alpha_1 \in \mathbb{F}_q \setminus \{0\}\) and~\(\alpha_i \in \mathbb{F}_q\) for~\(2 \leq i \leq m.\) Letting~\(Q_0 := Q(r_1,\ldots,r_m),\) we have \[\mathbb{P}(Q_0 =0) =  \sum_{j=2}^{m}\sum_{a_j \in \mathbb{F}_q} \mathbb{P}(Q_0 = 0 \mid r_2 = a_2,\ldots,r_m = a_m) \mathbb{P}(r_2 = a_2,\ldots,r_m = a_m)\] where
\[\mathbb{P}(Q_0 = 0 \mid r_2 = a_2,\ldots,r_m = a_m)  = \mathbb{P}\left(r_1 = -\frac{1}{\alpha_1} \left(\sum_{i=2}^{m} \alpha_i \cdot a_i\right)\right) = \frac{1}{q}.\] This implies that~\(\mathbb{P}(Q_0 = 0) = \frac{1}{q}\) and consequently we get
\begin{equation}\label{pm_est}
p_{1,m} =\frac{1}{q}
\end{equation}
proving~(\ref{pk_est_gen}) for~\(k=1.\)

To estimate~\(p_{k,m}\) for larger values of~\(k,\) we obtain a recursive relation for~\(p_{k,m}\)
as follows. Let~\(Q =Q(x_1,\ldots,x_m)\) be a homogenous polynomial of degree~\(k\) and write
\begin{equation}\label{q_eq}
Q = x_1 \cdot R + S,
\end{equation}
where~\(R = R(x_2,\ldots,x_m)\) is a homogenous polynomial of degree~\(k-1\) in the variables~\(x_2,\ldots,x_m\)
and~\(S = S(x_2,\ldots,x_m)\) is a homogenous polynomial of degree~\(k\) in the variables~\(x_2,\ldots,x_m.\)
If~\(Q_0 := Q(r_1,\ldots,r_m) = 0\) then either~\(R_0 := R(r_2,\ldots,r_m) = 0\) and~\(S_0 := S(r_2,\ldots,r_m) = 0\)
or~\(R_0 \neq 0\) and so we have
\begin{equation}\label{rec_eq}
\mathbb{P}(Q_0 = 0) = \mathbb{P}(R_0 = 0, S_0 = 0) + \mathbb{P}(R_0 \neq 0, Q_0 = 0).
\end{equation}

The term
\begin{equation}\label{gen_3}
\mathbb{P}(Q_0 = 0, R_0 \neq 0) = \sum_{a \in \mathbb{F}_q \setminus \{0\}} \sum_{b \in \mathbb{F}_q}\mathbb{P}(Q_0 = 0 \mid R_0 = a,S_0=b)\mathbb{P}(R_0 = a,S_0 = b),
\end{equation}
and so for any~\(a \in \mathbb{F}_q \setminus \{0\}, b \in \mathbb{F}_q\) we have from~(\ref{q_eq}) that
\begin{equation}\label{q_0_est}
\mathbb{P}(Q_0 = 0 \mid  R_0 = a,S_0=b) = \mathbb{P}\left(r_1 = -\frac{b}{a}\right) = \frac{1}{q}.
\end{equation}
Consequently
\begin{equation}\label{gen_zz}
\mathbb{P}(R_0 \neq 0, Q_0 = 0) = \frac{1}{q} \mathbb{P}(R_0 \neq 0).
\end{equation}

To obtain the lower bound in~(\ref{z_est}), we use~(\ref{rec_eq}) and~(\ref{gen_zz}) to get that
\begin{eqnarray}
\mathbb{P}(Q_0 = 0) &\leq& \mathbb{P}(R_0=0)+ \frac{1}{q} \mathbb{P}(R_0 \neq 0) \nonumber\\
&=& \frac{1}{q} + \left(1-\frac{1}{q}\right) \mathbb{P}(R_0 =0) \nonumber\\
&\leq& \frac{1}{q} + \left(1-\frac{1}{q}\right)p_{k-1,m-1}. \nonumber
\end{eqnarray}
From~(\ref{pk_def}) we therefore get that
\begin{equation}\label{gen_one2}
p_{k,m} \leq \frac{1}{q} + \left(1-\frac{1}{q}\right) p_{k-1,m-1}.
\end{equation}
Letting~\(\beta = 1-\frac{1}{q}\) and applying the recursion~(\ref{gen_one2}) repeatedly~\(i\) times, we get
\begin{eqnarray}
p_{k,m} &\leq& (1-\beta)(1+\beta +\beta^2 + \ldots \beta^{i}) + \beta^{i+1}\cdot p_{k-1-i,m-1-i} \nonumber\\
&=& (1-\beta^{i+1}) + \beta^{i+1} \cdot p_{k-1-i,m-1-i} \label{gen_three}
\end{eqnarray}
Setting~\(i = k-2\) in~(\ref{gen_three}) and using~(\ref{pm_est}) we then get~(\ref{pk_est_gen}).
This obtains the upper bound in~(\ref{z_est}).

For the lower bound in~(\ref{z_est}), we again use~(\ref{rec_eq}) and get
\begin{equation}\label{gne_uu}
\mathbb{P}(Q_0 = 0) \geq \mathbb{P}(R_0 \neq 0, Q_0 = 0) = \frac{1}{q}\mathbb{P}(R_0 \neq 0)
\end{equation}
by~(\ref{gen_zz}). Thus
\[\mathbb{P}(Q_0 = 0) \geq \frac{1}{q}(1-p_{k-1,m-1}) \geq \frac{1}{q} \left(1-\frac{1}{q}\right)^{k-1}\]
by~(\ref{pk_est_gen}) and this proves the lower bound in~(\ref{z_est}).~\(\qed\)

\setcounter{equation}{0}
\renewcommand\theequation{\thesection.\arabic{equation}}
\section{Applications of Theorem~\ref{thm1}} \label{examples}
\subsection*{Perfect matching}
In this subsection, we illustrate an application for the upper bound in~(\ref{z_est}). We recall the polynomial identity testing procedure to determine the presence of a perfect matching in a bipartite graph. Let~\(G = (X \cup Y,E)\) be a bipartite graph where~\(\#X = \#Y = k.\) A matching in~\(G\) is a set of edges that share no endpoint. A perfect matching in~\(G\) is a set of edges~\(M \subseteq E\) such that every vertex  is contained in exactly one edge of~\(M.\)

Let~\(A = [a_{i,j}]\) be a matrix whose rows are indexed by vertices in~\(X\) and whose columns are indexed by vertices in~\(Y\) with entries
\begin{equation}\label{a_entry}
a_{u,v} := \left\{
\begin{array}{cc}
x_{u,v}  & \text{ if } (u,v) \in E\\
0         & \text{ otherwise}
\end{array}
\right.
\end{equation}
where~\(\{x_{u,v}\}\) are distinct variables and~\((u,v)\) denotes the edge with endvertices~\(u\) and~\(v.\) It is well known (Theorem~7.3,~\cite{mot}) that the determinant of~\(A\) is zero if and only if the graph~\(G\) does not have a perfect matching.

Suppose~\(G\) has a perfect matching and we would like to use the above determinant criterion to devise a randomized algorithm for determining whether~\(G\) has a perfect matching or not. The Schwartz-Zippel procedure is as follows. Assuming that~\(\det(A)\) is a polynomial of degree~\(k,\) we first choose a field size~\(q \geq k+1.\) Fixing such a~\(q,\) we then  choose~\(\{x_{u,v}\}\) independent and identically distributed (i.i.d.) in~\(\mathbb{F}_q\) and compute the random determinant~\(det(A).\) If~\(det(A)=0\) we say that~\(G\) does not have a perfect matching else we say that~\(G\) has a perfect matching.

From the Schwartz-Zippel lemma we get that the probability that~\(det(A)\) is zero is at most~\(\frac{k}{q} < 1\) strictly and so the probability that our decision is wrong is at most~\(\frac{k}{q}.\) To reduce the probability of a wrong decision, we run the above procedure~\(n\) times using fresh independently random values for~\(\{x_{u,v}\}\) each time. If we get~\(det(A) =0\) all the~\(n\) times, we output~``\(G\) has no perfect matching" else we output~``\(G\) has a perfect matching". The probability that our decision is wrong in this case is at most~\(\left(\frac{k}{q}\right)^{n}\) which is small for all large~\(n,\) provided the field size~\(q \geq k+1.\)

Using Theorem~\ref{thm1} we now perform the above procedure using binary random variables. We first see that the determinant~\(det(A)\) of~\(A\) is a \emph{homogenous} polynomial of degree~\(k\) in the variables~\(\{x_{u,v}\}.\) We then set~\(x_{u,v}\) to be independent random binary values satisfying
\[\mathbb{P}(x_{u,v} = 0) = \frac{1}{2} = \mathbb{P}(x_{u,v} = 1)\] and evaluate~\(det(A).\) If~\(det(A) = 0,\) we output the statement~``\(G\) has no perfect matching"; else we output the statement~``\(G\) has a perfect matching". From~(\ref{z_est}), the probability that~\(det(A)\) equals zero is at most~\(1-\frac{1}{2^{k}}\) and so the probability that we output the wrong decision is at most~\(1-\frac{1}{2^{k}}.\)

As before, to reduce the probability of a wrong decision, we run the above procedure~\(n\) times using fresh independently random values for~\(\{x_{u,v}\}\) each time. If we get~\(det(A) =0\) all the~\(n\) times, we output~``\(G\) has no perfect matching" else we output~``\(G\) has a perfect matching". Again using~(\ref{z_est}), we get that the probability that our decision is wrong is at most~\(\left(1-\frac{1}{2^{k}}\right)^{n},\) which decays exponentially with~\(n.\)

Of course the tradeoff involved in the above procedure is the running time: our algorithm requires~\(n \times poly(k)\) running time since we need to compute~\(n\) determinants, each of size~\(k \times k.\) It would be interesting to design algorithms that require lesser computation and this would be a potential direction for future study.

\subsection*{Common zeros}
We illustrate the lower bound in~(\ref{z_est}) with an example involving common zeros of polynomials. We have the following result.
\begin{proposition}\label{prop1} Suppose~\(Q_1,Q_2,\ldots,Q_N\) are~\(N\) homogenous polynomials in the variables~\(x_1,\ldots,x_m,\) each with degree~\(k.\)  If~\(N \geq 1+q\cdot \left(\frac{q}{q-1}\right)^{k-1},\) then there are indices~\(1 \leq i \neq j \leq N\) such that~\(Q_i\) and~\(Q_j\) have a common zero.
\end{proposition}
We remark here that the Chevalley-Warning theorem~\cite{schmidt} is used to describe conditions under which a set polynomials whose sum degree is smaller than the number of variables, all have more than one common root. In Proposition~\ref{prop1}, we require that the number of polynomials is sufficiently large in order that at least two of the polynomials have a common root.

\emph{Proof of Proposition~\ref{prop1}}: We use the lower bound in~(\ref{z_est}) and prove by contradiction. Let~\({\cal Z}_i\) be the set of zeros of the polynomial~\(Q_i.\) If the sets~\(\{Z_i\}\) are all mutually disjoint, then by the lower bound in~(\ref{z_est}), the total number of elements in~\(\bigcup_{i=1}^{N} {\cal Z}_i\) is at least~\(N \cdot q^{m-k}(q-1)^{k-1} > q^{m}\) strictly, a contradiction.~\(\qed\)

\section*{\textit{Acknowledgment:}}

\textit{I thank Professors V. Arvind, C. R. Subramanian and the referee for crucial comments that led to an improvement of the paper. I also thank IMSc for my fellowships.}


\section*{\centering REFERENCES}
\setlength{\parindent}{0.32in}
\setlength{\leftskip}{0.5in}
\setlength{\parskip}{15pt}

\end{document}